# FUNDAMENTAL DOMAINS AND ANALYTIC CONTINUATION
# OF GENERAL DIRICHLET SERIES

**Dorin Ghisa**[1]

with an introduction by   **Les Ferry**[2]

**Abstract**: Fundamental domains are found for functions defined by general Dirichlet series and by using basic properties of conformal mappings the Great Riemann Hypothesis is studied.

## Introduction

My interest in the topic of this paper appeared when I became aware that John Derbyshare, with whom I was in the same class as an undergraduate, had written a book on the Riemann Hypothesis and that this had won plaudits  for its remarkable success in representing an abstruse topic to a lay audience. On discovering this area of mathematics for myself through different other sources, I found it utterly fascinating. More than a quarter of century before, I was taken a strong interest in the work of Peitgen, Richter, Douady and Hubbard concerning Julia, Fatou and Mandelbrot sets. From this I had learned the basic truth that, however obstruse might be the method of definition, the nature of analytic function is very tightly controlled in practice. I sensed intuitively at that time the Riemann Zeta function would have to be viewed in this light if a proof of Riemann Hypothesis was to be achieved.

Many mathematicians appear to require that any fundamental property of the Riemann Zeta function and in general of any L-function must be derived from the consideration of the Euler product and implicitly from some mysterious properties of prime numbers. My own view is that it is in fact the *functional equation* which expresses the key property with which any method of proof of the Riemann Hypothesis and its extensions must interact. It is through elementary handling of the functional equation that we perceive the need to derive the two non trivial zeros of the function having the same value of the imaginary part and straddling the critical line if Riemann Hypothesis is to be refuted. It is surely no coincidence that Dorin Ghisa's treatment of the Great Riemann Hypothesis (GRH) in this paper effectively examines axioms for the existence of functions satisfying a similar functional equation. Obviously, the functions in the Selberg class are valid members, yet he is considering functions obtained from more general Dirichlet series. This approach usefully distinguishes this paper from the more amateurish offerings which focus directly on the Riemann Hypothesis itself.

By the autumn of 2012, my interest in the topic had grown to an extent such that I had devised a crude method of my own for analyzing the behavior of the Riemann Zeta function $\zeta(\sigma + it)$ for large values of $t$. This incorporated a simple algorithm which would allow the search for potential counter examples to be dramatically narrowed. The exercise as a whole convinced me that the Riemann Hypothesis could be treated by using the considerable power of the methods based on *conformal mappings*, which abound in the whole field of complex dynamics. I realized that even if a proof of the Riemann Hypothesis were to be submitted and were to gain acceptance, my own acceptance of it would be lukewarm unless and until I could relate it to methods using properties of conformal mappings.

An appropriate Google search led immediately to the work of Dorin Ghisa. I realized that his work made careful use of basic properties of conformal mappings to establish a *natural partition* of the complex plane into a countably infinite collection of strips $S_k$ within which the behavior of the Riemann Zeta function and its derivative were carefully described. This was a kind of very efficient *divide and conquer* technique, since the behavior of the two functions in each strip $S_k$ was similar to the point that it was enough to study only one of these strips. Moreover, he found a just as *natural way* to divide the strips $S_k$ into fundamental domains, where $\zeta'(s) \neq 0$, hence the mapping realized by $\zeta(s)$ was bijective. Finally, the supposition that two non trivial zeros could exist in symmetric positions with respect to the critical line was refuted, as I had always expected it would need to be, by relatively simple topological arguments.

Since that time I have made contact with Professor Ghisa and he has kindly involved me in these matters. The method of generalization which he has adopted in this paper is I suspect applicable much more widely; he and I are looking at other possibilities.

Finally, I observe that Andrew Wiles' proof of Fermat Last Theorem represented a small and relatively insignificant part of the overall result which he was presenting at that time. Like wise the classic problem of the Riemann Hypothesis represents a particular case of the more far reaching result presented here.

Thus the worldwide body of academic mathematicians may be on the brink of a significant opportunity to declare a recognizable success to the wider community. If so, they should surely take it.

19-th January 2015

## 1. General Dirichlet Series

Let $\Lambda = \{0 = \lambda_1 \le \lambda_2 \le \ldots\}$ be an increasing sequence of non negative numbers such that $\lim_{n->\infty} \lambda_n = +\infty$ and let $A = \{a_n\}$ be an arbitrary sequence of complex numbers, $a_n \ne 0$ for infinitely many $n$. A series of the form

$$(1) \qquad \zeta_{A,\Lambda}(s) = \sum_{n=1}^{\infty} a_n e^{-\lambda_n s}$$

is called **general Dirichlet series**. We notice that when $\lambda_n = \ln n$, then $e^{-\lambda_n s} = \frac{1}{n^s}$ and (1) is an ordinary Dirichlet series and when $\lambda_n = n - 1$ we obtain a power series in $e^{-s}$. In fact, any power series

$$(2) \qquad \sum_{n=0}^{\infty} a_n (z - z_0)^n$$

can be converted into the general Dirichlet series

$$(3) \qquad \sum_{n=0}^{n} a_n e^{-ns}$$

by the substitution $z - z_0 = e^{-s}$. The Hadamard's formula

$$(4) \qquad 1/R = \limsup_{n->\infty} |a_n|^{1/n}$$

gives the radius of convergence of the series (2). It means that the series (2) converges absolutely for $|z - z_0| < R$ and diverges for $|z - z_0| > R$. Therefore the series (3) converges absolutely for $|e^{-\sigma - it}| < R$, i.e. for $\sigma > 1/R$ and diverges for $\sigma < 1/R$. The number $1/R$ is the *abscissa of convergence* of the series (3), which coincides with the abscissa of absolute convergence. For arbitrary Dirichlet series abscissa of convergence and that of absolute convergence can be different.

Namely (see [3], Theorem 8.2]), if (1) converges for some $s = \sigma + it$ with $\sigma > 0$, but diverges for all $s = \sigma + it$ with $\sigma < 0$, then with the notation

$$(5) \qquad \sigma_c = \limsup_{n->\infty} \ln|\sum_{k=1}^{n} a_k|^{1/\lambda_n} \ ,$$

$\zeta_{A,\Lambda}(s)$ converges for $\operatorname{Re} s > \sigma_c$ and diverges for $\operatorname{Re} s < \sigma_c$. The number $\sigma_c$ is called the abscissa of convergence of the series. Also, if

$$(6) \qquad \sigma_a = \limsup_{n->\infty} \ln[\sum_{k=1}^{n} |a_k|]^{1/\lambda_n}$$

then the number $\sigma_a$ is the *abscissa of absolute convergence* of the series. Obviously, $-\infty \le \sigma_c \le \sigma_a \le +\infty$.

It makes no sense to deal with general Dirichlet series having $\sigma_a = +\infty$. The condition $\sigma_a < +\infty$ appears to replace the Ramanujan condition (See [11]) for ordinary Dirichlet series.

When studying the zeros of $\zeta_{A,\Lambda}(s)$ we can suppose that $a_1 = 1$, since otherwise, if $a_1 = a_2 = \ldots = a_{m-1} = 0$, $a_m \notin \{0,1\}$ and $\lambda_m > 0$, we can study instead the series $\frac{e^{\lambda_m s}}{a_m}\zeta_{A,\Lambda}(s)$, which is still a general Dirichlet series and has the same abscissa of convergence and the same zeros as $\zeta_{A,\Lambda}$. For this last function we have

$$(7) \qquad \lim_{\sigma \to +\infty} \frac{e^{\lambda_m s}}{a_m}\zeta_{A,\Lambda}(\sigma + it) = 1$$

In the following we will study only functions $\zeta_{A,\Lambda}$ with $A = (1, a_2, \ldots)$, where $a_k \neq 0$ for infinitely many $k$ and for which

$$(8) \qquad \lim_{\sigma \to +\infty} \zeta_{A,\Lambda}(\sigma + it) = 1$$

It can be shown that $\zeta_{A,\Lambda}(\sigma + it)$ tends to 1 uniformly with respect to $t$ as $\sigma \to +\infty$. Indeed, for $\sigma \geq \sigma_0 > \sigma_a$ we have

$$|\Sigma_{n=2}^{\infty} a_n e^{-\lambda_n s}| \leq \Sigma_{n=2}^{\infty} |a_n| e^{-\lambda_n(\sigma - \sigma_0)} e^{-\lambda_n \sigma_0} = e^{-\lambda_2(\sigma - \sigma_0)} \Sigma_{n=2}^{\infty} |a_n| e^{-(\lambda_n - \lambda_2)(\sigma - \sigma_0)} e^{-\lambda_n \sigma_0} \leq Ce^{-\lambda_2 \sigma} \to 0 \text{ as } \sigma \to +\infty,$$

where $C$ does not depend on $\sigma$, given the fact that $e^{-(\lambda_n - \lambda_2)(\sigma - \sigma_0)} \leq 1$, $e^{\lambda_2 \sigma_0}$ is a constant and $\Sigma_{n=2}^{\infty} |a_n| e^{-\lambda_n \sigma_0}$ converges.

## 2. Analytic Continuation of General Dirichlet Series

The functions $\zeta_{A,\Lambda}(s)$ are analytic in the half plane $\operatorname{Re} s > \sigma_c$ due to the fact that the series (1) converges uniformly on compact subsets of this half plane. On the other hand, every term of the series has an essential singularity at $s = \infty$, hence the same is true for $\zeta_{A,\Lambda}(s)$. By Big Picard Theorem every complex value (except possibly a lacunary one) is taken infinitely many times by $\zeta_{A,\Lambda}(s)$ in any neighborhood of $s = \infty$. In particular, if $z = 0$ is not a lacunary value, the function $\zeta_{A,\Lambda}(s)$ has infinitely many zeros.

On the other hand, due to the fact that we have the limit (5) uniformly with respect to $t$, there is $\sigma_0 \geqslant \sigma_c$ such that $\zeta_{A,\Lambda}(\sigma + it) \neq 0$ for $\sigma \geq \sigma_0$. It is known that for the Dirichlet L-series there is no zero in the half plane of absolute convergence. However, there are infinitely many zeros in the strip $0 < \operatorname{Re} s < 1$ of conditional convergence (the so called non trivial zeros) and also countable many real zeros appear after the analytical continuation of the respective series to the whole complex plane (the trivial ones). The question arises if this is generally true, i.e. if it is true for arbitrary general Dirichlet series, or if not, then what are the series for which it is true.

Let us suppose that $\sigma_a$ defined in (6) for the series (1) is finite. We start from a point $s$ with $\operatorname{Re} s > \sigma_c$ and perform Weierstrass continuation of $\zeta_{A,\Lambda}(s)$. There are two possibilities:

a). Every point of the line $\operatorname{Re} s = \sigma_c$ is a singular point of the series (1) and therefore any Taylor series obtained by expanding the terms of (1) has a convergence disc strictly included in the half plane $\operatorname{Re} s > \sigma_c$ .

b). There are Taylor series obtained as in a) with the convergence disc overlapping the half plane $\operatorname{Re} s < \sigma_c$.

Obviously, in the case a) the analytic continuation of the series across the line $\operatorname{Re} s = \sigma_c$ is impossible. An example of such a series can be obtained by converting the Hadamard's series $1 + \Sigma_{n=0}^{\infty} z^{2^n}$ into a general Dirichlet series, namely $1 + \Sigma_{n=0}^{\infty} e^{-2^n s}$. It is known that every point of the circle $|z| = 1$ is a singular point for the Hadamard's series, therefore every point of the line $\operatorname{Re} s = 0$ is a singular point for the corresponding Dirichlet series.

There are no known conditions on the sequences $A$ and $\Lambda$ allowing one to decide if the series (1) belongs to the case a) or b) above. If it does belong to the case b), then two situations can arise:

(i). The continuation takes place into the whole half plane $\operatorname{Re} s \leq \sigma_c$, except possibly for a discrete set of poles, giving rise to a meromorphic function which is locally injective, except for a discrete set of points, the so called branch points. We keep the notation $\zeta_{A,\Lambda}(s)$ for this extended function.

(ii). The singular points in the half plane $\operatorname{Re} s \leq \sigma_c$ form a continuum and/or some essential singular points exist.

An example of the case (ii) function can be obtained by converting a certain Blaschke product into a Dirichlet series.

Take as zeros of the respective Blaschke product the numbers of the form

$a_{n,k} = (1 - \frac{1}{3^n})e^{2k\pi i/2^n}, \; k = 1, 2, \ldots, 2^n.$

Then $\Sigma_{n,k}(1 - |a_{n,k}|) = 2$ and by Blaschke criterion the respective Blaschke product is convergent in the unit disc. However, every point of the unit circle is a singular point for the sum of this series, since it is a limit point of the poles $1/\overline{a}_{n,k}$.

By a formal computation we get:

(9)     $\Pi_{n,k}\frac{z-a_{n,k}}{1-\overline{a}_{n,k}z} = \Pi_{n,k}(z-a_{n,k})(1+\overline{a}_{n,k}z+\ldots) = \Pi_{n,k}[-a_{n,k} + (1-|a_{n,k}|^2)z+\ldots] = \Sigma_{n=0}^{\infty}\alpha_n z^n =$

$1 + \Sigma_{n=1}^{\infty}\alpha_n e^{-ns}$, where $z = e^{-s}$. We can choose $k$ such that the limit points of poles form a continuum, hence the singular points of the Dirichlet series are like in (ii).

We will deal in the following only with the case (i) and some of the theorems will adopt assumptions specific to the Selberg class of *standard* Dirichlet series ([7], [10], [11], [14]). In particular, it will be convenient, although not always necessary, to postulate that $\sigma_a = 1$ and the existence of at most one simple pole at $s = 1$. It is useful to consider the couple $(\zeta_{A,\Lambda}, \mathbb{C})$ as a *ramified (or branched) covering Riemann surface of* $\overline{\mathbb{C}}$ (See [2], p. 31). *Continuations along curves, or lifting of curves*, as defined in that book, will play in what follows an important role. When in the process of continuation along a curve we meet a pole, the continuation stops, meaning that the image of a finite arc starting at the pole is an unbounded curve and when transcending the pole, a new curve is obtained. This appears as a kind of embarrassment which can be avoided by the respective assumption. When a branch point is met, the continuation can follow from there on a finite number of different paths such that every one of them is mapped bijectively by $\zeta_{A,\Lambda}(s)$ onto the same curve. We will call a *component* of the pre-image of a curve $\gamma$ any curve obtained by continuation along the whole curve $\gamma$ starting from a point $s$ with $\zeta_{A,\Lambda}(s) \in \gamma$.

Thus, when branch points are met we are forced to allow some components of the pre-image of a curve to have common parts. This is another embarrassment to be avoided, when possible. It is known that a branch point $s_0$ of multiplicity $q \geq 2$ of $\zeta_{A,\Lambda}(s)$ is either a zero of order $q - 1$ of $\zeta'_{A,\Lambda}(s)$ or a pole of order $q$ of $\zeta_{A,\Lambda}(s)$. In a neighborhood of $s_0$ we have $\zeta_{A,\Lambda}(s) = (s - s_0)^q \varphi(s)$, in the first case and $\zeta_{A,\Lambda}(s) = (s - s_0)^{-q} \varphi(s)$ in the second case, where $\varphi$ is analytic at $s_0$ and $\varphi(s_0) \neq 0$ (See [1], p.133 and [9], p.7).

The pre-image of any curve $\gamma$ passing through $\zeta_{A,\Lambda}(s_0)$ $(\zeta_{A,\Lambda}(s) = \infty$ if $s_0$ is a pole) is formed with $q$ curves passing through $s_0$ two consecutive of which make at $s_0$ an angle of $\pi/q$. We will perform continuations along two types of curves: circles centered at the origin and rays issuing from the origin. They form an orthogonal net which is the image by $\zeta_{A,\Lambda}(s)$ extended to the whole plane of a net which is orthogonal, except at the branch points of this function. If a branch point is of order two, then the components issuing from that point are still orthogonal two by two, as seen in Fig. 1 below.

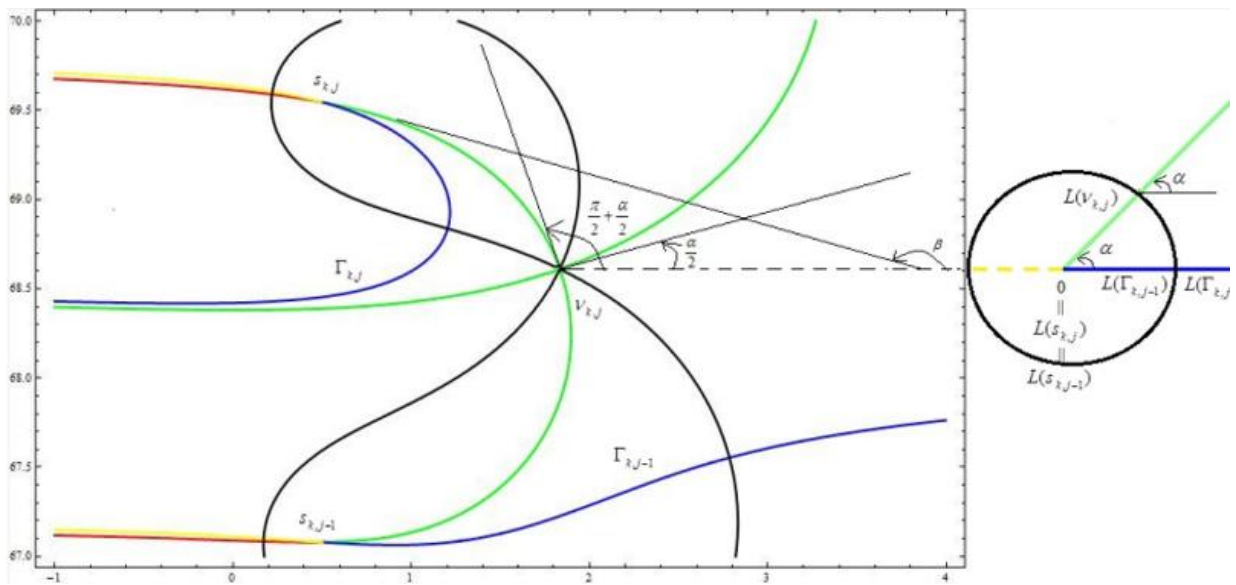

Fig. 1

Due to the continuity of $\zeta_{A,\Lambda}(s)$ the pre-image by this function of a small open disc centered at the origin is a collection of connected open sets containing each one a unique zero of $\zeta_{A,\Lambda}(s)$. Since the branch points form a discrete set, we can choose the radius of the disc in such a way that the pre-image of the respective circle does not contain any branch point. Then the boundary of such a set is a curve orthogonal to the components of the pre-image of any ray issuing from the origin.

When increasing the radius of the respective small disc, those open sets expand until the boundaries of two or more of them touch at a point $v$. This is a branch point of $\zeta_{A,\Lambda}(s)$, since the injectivity at $v$ is violated, hence $\zeta'_{A,\Lambda}(v) = 0$. Indeed, we can see that $v \neq 1$, since $s = 1$ is either a regular point or a simple pole and $\zeta_{A,\Lambda}(s)$ is injective at that point. Increasing more that radius the respective components will fuse into a unique connected component and the process continues. Unbounded components can be obtained in this way.

**Theorem 1**. *Suppose that the series (1) can be extended to a meromorphic function in the whole complex plane. We keep the notation $\zeta_{A,\Lambda}(s)$ for this function. Then the pre-image of the closed unit disc by $\zeta_{A,\Lambda}(s)$ has at least one unbounded connected component.*

*Proof:* Assume that all the connected components of the pre-image of the closed unit disc were bounded. Obviously, this pre-image cannot contain the pole of the function, since its image is a bounded set, therefore $\zeta_{A,\Lambda}(s)$ is continuous on the respective pre-image. Due to the fact that the limit (8) is taken uniformly with respect to $t$, for every $\epsilon > 0$ there is $\sigma_\epsilon$ such that $\sigma > \sigma_\epsilon$ implies $|\zeta_{A,\Lambda}(\sigma + it) - 1| < \epsilon$. Hence the half-plane $U_\epsilon = \{\sigma + it \mid \sigma > \sigma_\epsilon\}$ is mapped by $\zeta_{A,\Lambda}(s)$ into the disc centered at $z = 1$ and of radius $\epsilon$. For $\delta > 0$, let us take a $\delta$-neighborhood $V$ of the pre-image of the closed unit disc, whose connected components are bounded open sets. Due to the fact that $\zeta_{A,\Lambda}(s)$ is an open mapping, $V$ is mapped by $\zeta_{A,\Lambda}(s)$ onto an open set containing the closed unit disc. If $\epsilon > 0$ is small enough, the closed disc centered at the origin and of radius $1 + \epsilon$ is included in that open set hence its pre-image is included in $V$, therefore it has only bounded connected components. However, one of them should contain the unbounded set $U_\epsilon$ and this is a contradiction.

Any connected component of the pre-image of the unit disc contains at least one zero of the function. We will see later that the number of zeros in every such component is finite. Having in view the relation (8), as $z$ approaches 1 on the unit circle, for any point $s$ on the boundary of an unbounded component $\Delta_k$ of its pre-image such that $\zeta_{A,\Lambda}(s) = z$ we must have $\mathrm{Re}\, s \to +\infty$. Thus, that boundary has the shape of a parabola with the branches extending to infinity in the right half plane as $z \to 1$ on the unit circle from the left and from the right. We will show later that in fact the pre-image by $\zeta_{A,\Lambda}(s)$ of the unit disc has infinitely many unbounded connected components. Four of them can be seen in Fig. 2 below colored part red, part white. The curve separating the two colors is mapped bijectively by $\zeta_{A,\Lambda}(s)$ onto the diameter $(-1, 1)$. The interval $(1, +\infty)$ must be in turn the image of some curve $\Gamma'_k$ exterior to $\Delta_{kk}$.

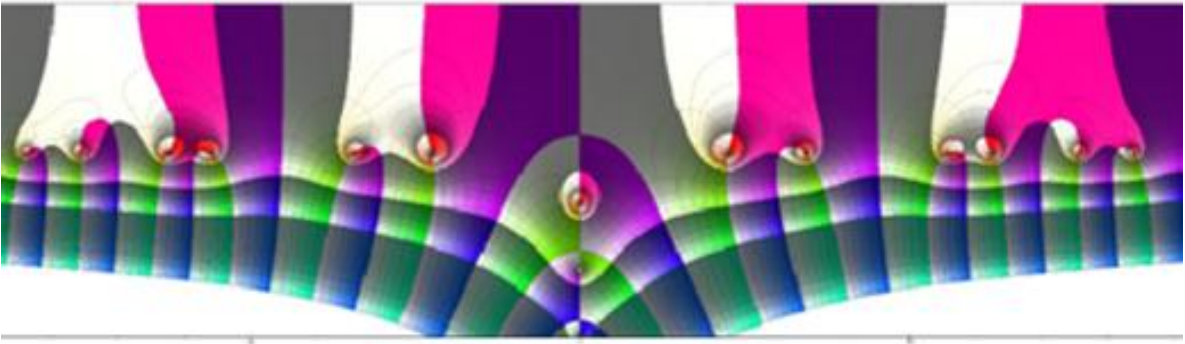

Fig.2

**Theorem 2**. *Every unbounded component $\Delta_k$ of the pre-image of the unit disc is situated between two consecutive components $\Gamma_k'$ and $\Gamma_{k+1}'$ of the pre-image of $\mathbb{R}$ which are mapped bijectively onto the interval $(1, +\infty)$ and vice-versa, there is a unique unbounded connected component of the pre-image of the unit disc between two consecutive such curves. That set contains a unique unbounded component $\gamma_{k,0}$ of the pre-image of the interval $(0, 1)$ such that for $\sigma + it \in \gamma_{k,0}$ we have $\lim_{\sigma \to +\infty} \zeta_{A,\Lambda}(\sigma + it) = 1$.*

*Proof:* Let us deal first with the pre-image of a ray $\eta_\alpha$ starting from the origin and making a small angle $\alpha$ with the positive real half axis. If $z = \zeta_{A,\Lambda}(s) \in \eta_\alpha$ then we cannot have $\operatorname{Re} s \to +\infty$, since $z$ cannot tend to 1 on $\eta_\alpha$. On the other hand, the continuation along $\eta_\alpha$ starting from any zero of $\zeta_{A,\Lambda}(s)$ is unlimited if $\eta_\alpha$ does not meet the pole and it is unique if $\eta_\alpha$ does not meet any branch point of $\zeta_{A,\Lambda}(s)$. The component $\Delta_k$ cannot contain a pole since its image is bounded. On the other hand, since the branch points form a discrete set, $\eta_\alpha$ can avoid their image, therefore the continuation along $\eta_\alpha$ from any zero of $\zeta_{A,\Lambda}(s)$ is unlimited and unique. We denote by $\Gamma_\alpha$ a component of the pre-image of $\eta_\alpha$ intersecting the boundary of $\Delta_k$ at a point $s_\alpha$ such that $\zeta_{A,\Lambda}(s_\alpha) = e^{i\alpha}$ and by $s_0$ the corresponding zero of $\zeta_{A,\Lambda}(s)$. Due to the continuity of $\zeta_{A,\Lambda}(s)$, letting $\alpha \to 0$ through positive and negative values the curve $\Gamma_\alpha$ must approach components of the pre-image of the intervals $(0, 1)$ and $(1, +\infty)$, the first situated inside $\Delta_k$ and the others outside $\Delta_k$. By the monodromy theorem there is a unique component $\gamma_{k,0}$ of the pre-image of the interval $(0, 1)$ obtained by continuation along this interval starting from $s_0$ and $\Gamma_\alpha$ approaches this curve as $\alpha$ tends to zero by positive or negative values. On the other hand, the part of $\Gamma_\alpha$ situated outside $\Delta_k$ is on one side or on the other side of $\Delta_k$ depending whether $\alpha$ is positive or negative. Therefore, we get two components $\Gamma_k'$ and $\Gamma_{k+1}'$ of the pre-image of the interval $(1, +\infty)$ as the limit positions of $\Gamma_\alpha$ as $\alpha$ tends to zero by positive, respectively negative values. On both these components we have $\lim_{\sigma \to +\infty} \zeta_{A,\Lambda}(\sigma + it) = 1$. The continuation from $s_0$ along $\eta_\alpha$ meets a unique point $s_\alpha$ such that $\zeta_{A,\Lambda}(s_\alpha) = e^{i\alpha}$, thus there can be just one component $\Delta_k$ with this property.

*Note*: The continuation from $s_0$ along the negative real half axis is also unlimited, as long as it doesn't meet the pole, so with this exception, we obtain a unique component $\Gamma_{k,0}$ of the pre-image of the real axis which is projected by $\zeta_{A,\Lambda}(s)$ onto the interval $(-\infty, 1)$ such that $\lim_{\sigma \to +\infty} \zeta_{A,\Lambda}(\sigma + it) = 1$ for $\sigma + it \in \Gamma_{k,0}$ and two curves $\Gamma'_k$ and $\Gamma'_{k+1}$ on which $\lim_{\sigma \to +\infty} \zeta_{A,\Lambda}(\sigma + it) = 1$. We will see later that on each one of these three curves $\sigma$ can take any negative value. Fig. 3 below illustrates this situation for the particular case of the Riemann Zeta function which has been implemented on *Mathematica* software. Some computation has been done for more general L-functions (see [13]).

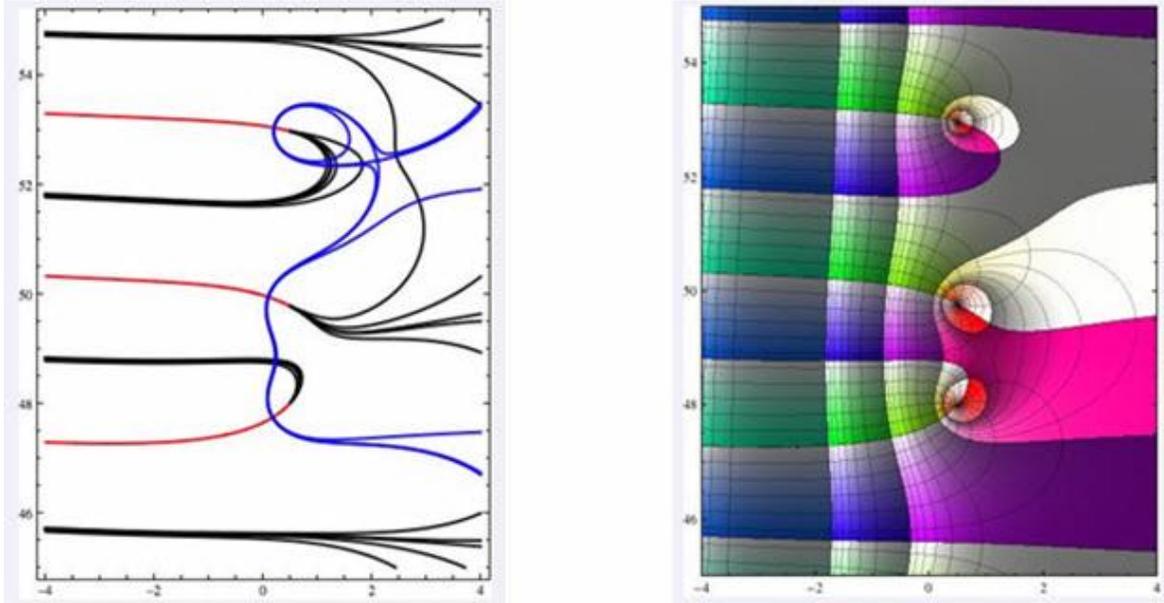

Fig.3

**Theorem 3**. *No couple of curves $\Gamma'_k$ and $\Gamma'_{k+1}$, with one exception, can meet each other. Thus, they form infinite strips $S_k$ containing at least one zero of $\zeta_{A,\Lambda}(s)$.*

*Proof:* Let us deal first with the possible exception. This is the case when the continuation along the interval $(1, +\infty)$ hits the pole $s = 1$. It happens in the case illustrated by Fig. 2 in which a part of the real axis belongs to $\Gamma'_0$ and $\Gamma'_1$. Other curves $\Gamma'_k$ and $\Gamma'_{k+1}$ cannot meet each other. Indeed, if they met at a point $s_0$, then they would bound a domain which is mapped conformally by $\zeta_{A,\Lambda}(s)$ onto the whole complex plane with a slit alongside the segment $[1, \zeta_{A,\Lambda}(s_0)]$ if $\zeta_{A,\Lambda}(s_0) > 1$, respectively $[\zeta_{A,\Lambda}(s), 1]$ if $\zeta_{A,\Lambda}(s_0) < 1$. That domain should contain the pole, which is impossible.

Obviously, no $\Gamma'_k$ can meet the pre-image of the interval $(-\infty, 1)$ since $(-\infty, 1) \cap (1, +\infty) = \varnothing$. We allow $k$ to take any integer value such that $\Gamma'_k$ is below $\Gamma'_{k+1}$ for every $k$ and $S_0$ is the strip containing the point $s = 1$. For the Riemann Zeta function $S_0$ contains all the trivial zeros and therefore infinitely many components of the pre-image of the real axis. On the other hand, a Dirichlet L-function defined by an imprimitive Dirichlet character has infinitely many imaginary trivial zeros (see [4] and [8]) which therefore do not belong to $S_0$. The geometry of the pre-image of the real axis is illustrated in Fig. 4 below for the case of two Dirichlet L-functions the first of a complex Dirichlet character and the second of a real one.

**Theorem 4**. *There are infinitely many strips $S_k$ covering the whole complex plane.*

*Proof:* Suppose that above a curve $\Gamma'_{k_0}$ there is no other curve $\Gamma'_k$. Let $s_0$ be a point above $\Gamma'_{k_0}$ and let $z_0 = \zeta_{A,\Lambda}(s_0)$. For any given point $z \in \mathbb{C}$ there is a curve $\eta$ having in common with the interval $[1, \infty)$ at most the point $z$ and which connects $z_0$ and $z$. The continuation along $\eta$ from $s_0$ brings us to a point $s$ such that $\zeta_{A,\Lambda}(s) = z$. Thus the closed domain above $\Gamma_{k_0}$ is mapped by $\zeta_{A,\Lambda}(s)$ locally injectively, except possibly for a discrete set of points, onto the whole complex plane with $\Gamma'_{k_0}$ mapped bijectively onto the interval $[1, +\infty)$. Let $x_0 \in (1, +\infty)$. Then there is $s_0 \in \Gamma'_{k_0}$ such that $\zeta_{A,\Lambda}(s_0) = x_0$. We can choose $x_0$ such that $\zeta'_{A,\Lambda}(s_0) \neq 0$. However, the local inverse of $\zeta_{A,\Lambda}(s)$ cannot be injective at $x_0$, which is a contradiction. Consequently, there should be infinitely many curves $\Gamma'_k$ above any curve $\Gamma'_{k_0}$ and, by symmetry, below it.

Suppose that we cannot have $\sigma \to -\infty$ on a $\Gamma'_k$, i.e. there is $\sigma_0$ such that $\sigma + it \in \Gamma'_k$ implies $\sigma > \sigma_0$. Then the same is true for any $\Gamma'_j$, $j \neq k$. Indeed, there is a zero of $\zeta_{A,\Lambda}(s)$ between $\Gamma'_k$ and $\Gamma'_j$ and therefore a component $\gamma$ of the pre-image of the interval $(0,1)$. Let $s_1 = \sigma_1 + it_1$ be such that $\sigma_1 < \sigma_0$ and let $s_2 \in \gamma$. Connect $\zeta_{A,\Lambda}(s_1)$ and $\zeta_{A,\Lambda}(s_2)$ by an arc $\eta$ avoiding the interval $(1, +\infty)$. Then the continuation along $\eta$ from $s_1$ should avoid both $\Gamma'_k$ and $\Gamma'_j$, which is impossible. Consequently we can let $\sigma \to -\infty$ on every $\Gamma'_k$ and on every $\Gamma_{k,j}$ as stated in the *Note* above and seen in Fig.4 below.

We have on all these curves

$$(10) \qquad \lim_{\sigma \to -\infty} \zeta_{A,\Lambda}(\sigma + it) = \pm\infty$$

We notice that such a configuration should be obtained for any general Dirichlet series (1) for which $\sigma_c$ given by (5) is finite, the relation (8) is true and which can be continued to a meromorphic function in the whole plane with at most one pole at $s = 1$. For the next theorems we need to do another assumption which is true for the Selberg class, namely that $\zeta_{A,\Lambda}(s)$ verifies a functional equation. Those equations take particular forms for particular L-functions all connecting the values of the respective function at $s$ and $1 - s$ by means of a multiplier. If we label the zeros of the multiplier as *trivial* zeros, and the remaining ones as *non trivial,* then the functional equation guarantees that the non trivial zeros appear in couples of the form $\sigma + it$ and $1 - \sigma + it$. The Grand Riemann Hypothesis says that for any such couple we have necessarily $\sigma = 1 - \sigma$, i.e. $\sigma = 1/2$. Thus we postulate that $\zeta_{A,\Lambda}(s)$ satisfies such a functional equation, without specifying any particular form of the multiplier.

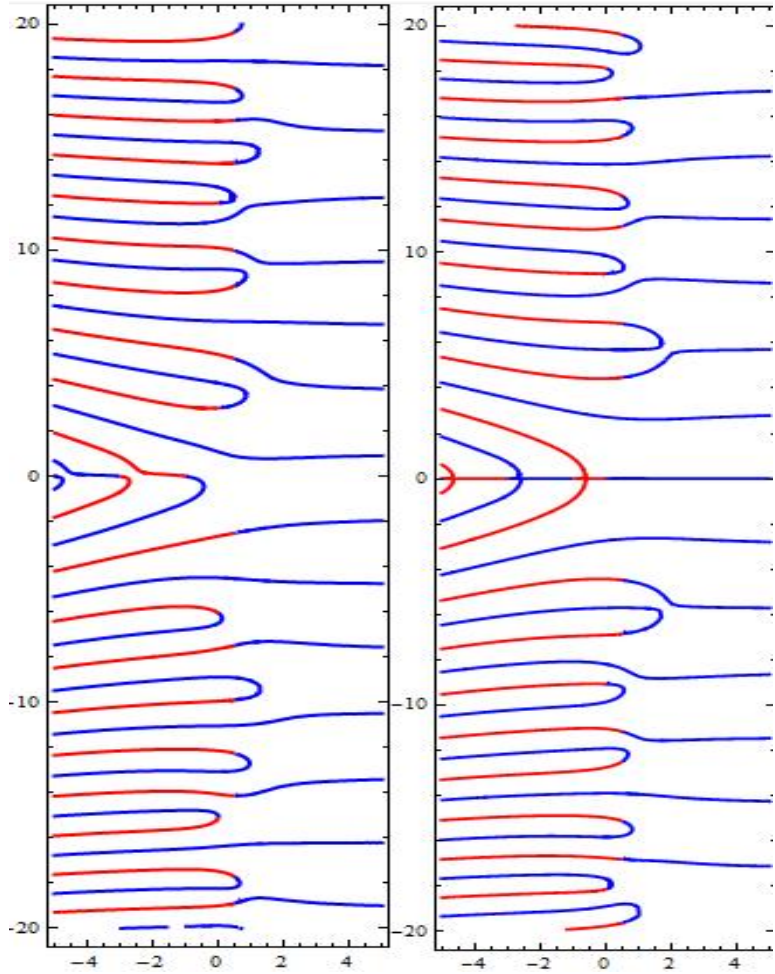

Fig.4

**Theorem 5**. *Every strip $S_k$, $k \neq 0$ contains a finite number $j_k$ of non trivial zeros of $\zeta_{A,\Lambda}(s)$ and $j_k - 1$ zeros of $\zeta'_{A,\Lambda}(s)$.*

*Proof:* Indeed, due to the relation (8), there is $\sigma_0 > 1$ such that $\zeta_{A,\Lambda}(\sigma + it) \neq 0$ for $\sigma > \sigma_0$. Due to the functional equation, the non trivial zeros are all in a strip $\{\sigma + it \parallel 1 - \sigma_0 < \sigma < \sigma_0\}$. The intersection of this strip with $S_k$ is a bounded set which can contain only a finite number $j_k$ of zeros. In particular, the number of zeros belonging to any connected component of the pre-image of the unit disc should be finite, as previously stated. On the other hand, for a given $k \neq 0$, we can take $r > 0$ small enough such that the components of the pre-image of the disc centered at the origin and of radius $r$ situated in $S_k$ do not overlap. Letting $r$ increase, these components expand, such that for an $r = r_0$ two of them touch each other at a point $v$. This is a branch point of $\zeta_{A,\Lambda}(s)$, hence $\zeta'_{A,\Lambda}(v) = 0$. Continuing to increase $r$ some other branch points can be obtained. In fact, all the zeros of $\zeta'_{A,\Lambda}(s)$ from $S_k$ are obtained in this way and a complete binary tree can be formed having as leafs the zeros of $\zeta_{A,\Lambda}(s)$ and as internal nodes the zeros of $\zeta'_{A,\Lambda}(s)$. When increasing $r$ past 1 the points of the pre-image of $r$ are all on the pre-image of the interval $(1, +\infty)$, hence turning indefinitely around the origin on $|z| = r$ will generate as component of the pre-image an unbounded curve intersecting all the curves $\Gamma'_k$. The contact of this curve with another component of the pre-image of the circle $|z| = r$ represents one more internal node of the respective complete binary tree. It is known that the number of these last ones must be $j_k - 1$. We can see this phenomenon in Fig. 5 below, where for the Riemann Zeta function the strip $S_1$ is illustrated with its unbounded component of the pre-image of the unit circle containing two non trivial zeros of the function and one zero of the derivative.

For the moment, we have to assume the possibility of multiple zeros and the previous numbers are obtained counting multiplicities. However, we will see soon that all of these zeros are simple. Let us notice that one zero of $\zeta_{A,\Lambda}(s)$ is situated on the unique component $\Gamma_{k,0}$ of the pre-image of the real axis belonging to $S_k$ which is mapped bijectively by $\zeta_{A,\Lambda}(s)$ onto the interval $(-\infty, 1)$, while the other $j_k - 1$ zeros from $S_k$ are situated on components $\Gamma_{k,j}$, $j \neq 0$ which are mapped bijectively by $\zeta_{A,\Lambda}(s)$ onto the whole real axis. If we color differently the pre-image of  the negative and that of the positive  real half axis, the zeros are at the junction of the two colors. When a point moves in the same direction on a small circle centered at the origin, any point from its pre-image will move around a zero of the function meeting alternately the two colors. We call this simple topological fact the *color alternating rule.* If we deal now with the pre-image of the real axis by $\zeta'_{A,\Lambda}(s)$, we will find a similar configuration with curves $\Upsilon'_k$ forming infinite strips and curves $\Upsilon_{k,j}$, $j \neq 0$ containing the zeros of the function, except that  the curves $\Upsilon_{k,0}$ do not contain any zero, since $\lim_{\sigma \to +\infty} \zeta'_{A,\Lambda}(s)$ is 0 and not 1. The effect of this is that every $\Upsilon_{k,0}$ is mapped bijectively by $\zeta'_{A,\Lambda}(s)$ onto the positive real half axis, every $\Upsilon'_k$ is mapped onto the negative real half axis and every $\Upsilon_{k,j}$, $j \neq 0$ is mapped bijectively  onto the whole real axis. If we use four different colors denoted **a**, **b**, **c**, **d** for the pre-image by $\zeta_{A,\Lambda}(s)$ and by $\zeta'_{A,\Lambda}(s)$ of the negative and respectively positive real half axis, we realize  that the intertwining curves must have specific colors (see [9], p. 102) and this is another simple topological fact which has been called the *color matching rule*. Namely, the color **b** meets always **c** and if $j \neq 0$, or $\sigma < 1/2$, the color **a** meets always **d**. Notice that for $\sigma > 1/2$ a curve $\Gamma_{k,0}$ can intersect $\Upsilon_{k,0}$ and therefore the color **b** will meet **d**. Such an exception has no bearing on the theorems which follow. A corollary of the two rules is the fact that the zeros of $\zeta_{A,\Lambda}(s)$ and those of $\zeta'_{A,\Lambda}(s)$ are all simple zeros. The proof of this affirmation for an arbitrary general Dirichlet series verifying the hypothesis of Theorems 1 and 2 is similar to that presented for the Riemann Zeta function in [9] and we will omit it.

It has been  shown in [4] that if a zero $s_k$ of the Dirichlet L-function $L(s ; \chi)$ generates in this way a zero $v_j$ of $L'(s ; \chi)$ then $\mathrm{Re}\, s_k < \mathrm{Re}\, v_j$. The proof uses again only facts mentioned in the Theorems 1 and 2 and  therefore the affirmation is true for general Dirichlet series. An illustration of this fact can be seen  in Fig. 1.

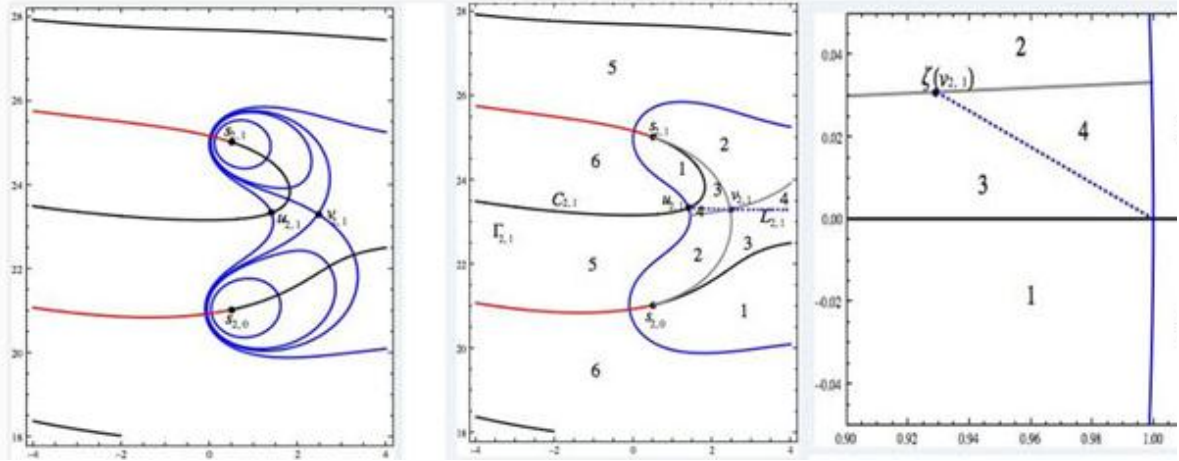

Fig. 5

## 3. Fundamental Domains of $\zeta_{A,\Lambda}(s)$

Every strip $S_k$, $k \neq 0$ is mapped by $\zeta_{A,\Lambda}(s)$ onto the whole complex plane with a slit alongside the interval $[1, +\infty)$ of the real axis. The mapping is $j_k$ to one, where $j_k$ is the number of zeros in $S_k$ of $\zeta_{A,\Lambda}(s)$. Infinitely many points from $S_0$ have the same image. A simple construction allows us to partition every $S_k$ into sub-strips which are mapped conformally by the function onto the whole complex plane with a slit. For $k \neq 0$, it consists in taking the pre-image of the segment between $z = 1$ and $\zeta_{A,\Lambda}(v_{k,j})$, where $v_{k,j}$ are the zeros of $\zeta'_{A,\Lambda}(s)$ situated in $S_k$. We formulate this result as a theorem whose proof is similar to that given in [9] and there is no need to repeat it.

**Theorem 6**. If $\zeta'_{A,\Lambda}(v_{k,j}) = 0$, then the pre-image of the segments $\eta_{k,j}$ from $z = 1$ to $z = \zeta_{A,\Lambda}(v_{k,j})$ has $j_{k-1}$ components situated in $S_k$ bounding together with the components of the pre-image of the interval $(1, +\infty)$ exactly $j_k$ fundamental domains $\Omega_{k,j}$ of $\zeta_{A,\Lambda}(s)$. Moreover, $\zeta_{A,\Lambda}(s)$ maps conformally every $\Omega_{k,j}$ onto the complex plane with the slit $\eta_{k,j} \cup (1, +\infty)$.

This can be seen in Fig. 5 above which illustrates the fundamental domains from $S_2$ of the Riemann Zeta function and the way they are mapped conformally onto the complex plane with the slit $\eta_{2,1} \cup (1, +\infty)$.

**Theorem 7**. The components of the pre-images of the real axis by $\zeta_{A,\Lambda}(s)$ and $\zeta'_{A,\Lambda}(s)$ are associated in pairs which intersect each other in points in which the tangents to the first are horizontal.

Proof: Let $\Gamma$ be a component of the pre-image by $\zeta_{A,\Lambda}(s)$ of the real axis. Let $s = s(x)$ be its parametric equation such that $\zeta_{A,\Lambda}(s(x)) = x$. Then $\zeta'_{A,\Lambda}(s(x))s'(x) = 1$, thus $\arg\zeta'_{A,\Lambda}(s(x)) + \arg s'(x) = 0 \pmod{2\pi}$. This equality implies that $\arg s'(x) = 0$ if and only if $\arg\zeta'_{A,\Lambda}(s(x)) = 0$ and $\arg s'(x) = \pi$ if and only if $\arg\zeta'_{A,\Lambda}(s(x)) = \pi$. By $\arg z$ we understand here the angle between the positive real half axis and the ray from $0$ to $z$. It is a simple geometric fact that $\Gamma$ must have at least one point $s(x)$ on which $\arg s'(x) = 0 \pmod{\pi}$. The affirmation becomes obvious if we look at the stereographic projections of $\Gamma$ and of its horizontal tangents onto the Riemann sphere.

At such a point, color **b** meets color **c** or color **a** meets color **d** or $\Gamma$ is $\Gamma_{k,0}$ and $\mathrm{Re}\, s(x) > 1/2$, in which case the intertwining curve is $\Upsilon_{k,0}$ and color **b** can meet color **d**. Due to the color alternating rule, no other curve $\Upsilon$ can intersect $\Gamma$, since at such an intersection point the color matching rule would be violated. Fig. 6 below illustrates these affirmations.

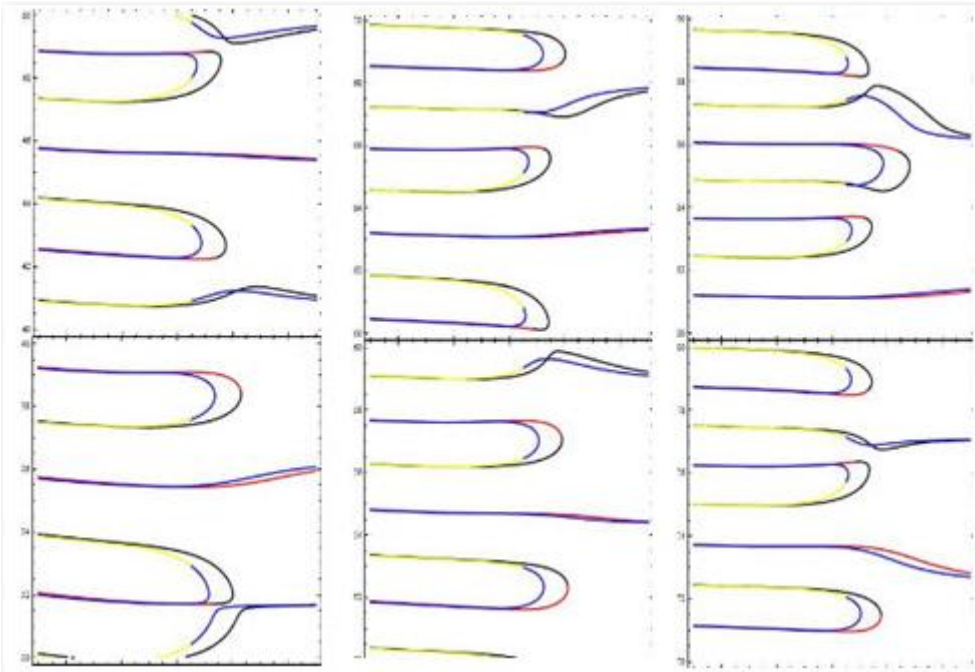

Fig. 6

## 4. The Grand Riemann Hypothesis

The existence of the continuation of a Dirichlet L-series to the whole complex plane and the existence of a functional equation are not independent facts. The same is true for more general L-series (see [12]) appearing in the adelic setting initiated by Tate and Weil and for which the **Grand Riemann Hypothesis** (**GRH**) has been formulated, namely that all the non trivial zeros of such an L-series lie on the line $\mathrm{Re}\, s = 1/2$. On the other hand, as noticed by [5], "there are L-functions, such as those attached to Maass waveforms, which do not seem to originate from geometry, and for which we still expect a Riemann hypothesis to be valid". In [6] the expectation is expressed that "if a functional equation and Euler product exists, then it is likely a type of Riemann hypothesis will hold". If we would like to reduce at minimum the requirements for a general Dirichlet series such that the GRH be valid, we can even get rid of stipulations regarding the Euler products, unless they do not appear implicitly in the functional equation. What we are proving now might be more than any form of GRH formulated up to date and something which could be generalized even more if the hypothesis of theorems 1 and 2 are relaxed .

**Theorem 8**. *Assume that the series (1) can be extended analytically to the whole complex plane, except for one possible pole at* $s = 1$, *that the extended function verifies the relation (8) and a functional equation of the form* $\zeta_{A,\Lambda}(s) = M(s)\,\overline{\zeta_{A,\Lambda}(1-\bar{s})}$. *Then the zeros of* $\zeta_{A,\Lambda}(s)$ *which are not zeros of* $M(s)$ *lie all on the line* $\operatorname{Re} s = 1/2$.

*Proof:* We notice that if $M(s_0) = 0$ then $\zeta_{A,\Lambda}(s_0) = 0$. We call $s_0$ a trivial zero of $\zeta_{A,\Lambda}(s)$. On the other hand, if $\zeta_{A,\Lambda}(1 - \sigma_0 + it_0) = 0$ then we should have also $\zeta_{A,\Lambda}(\sigma_0 + it_0) = 0$. We call $\sigma_0 + it_0$, respectively $1 - \sigma_0 + it_0$ non trivial zeros of $\zeta_{A,\Lambda}(s)$. The theorem states that for such a couple of zeros we have necessarily $\sigma_0 = 1 - \sigma_0$, since $\sigma_0 = 1/2$. In order to prove GRH for $\zeta_{A,\Lambda}(s)$, it is enough to show that if $s_1 = \sigma + it$ and $s_2 = 1 - \sigma + it$ are two non trivial zeros of $\zeta_{A,\Lambda}(s)$ then $s_1 = s_2$, therefore $\sigma = 1/2$. Suppose that two distinct zeros of this form exist and let $I$ be the segment connecting them. The parametric equation of $I$ is $s(\lambda) = (1 - \lambda)s_1 + \lambda s_2$, $0 \leq \lambda \leq 1$. Let $\gamma$ and $\gamma'$ be the images of $I$ by $\zeta_{A,\Lambda}(s)$, respectively by $\zeta'_{A,\Lambda}(s)$. Their parametric equations are respectively: $z(\lambda) = \zeta_{A,\Lambda}((1 - \lambda)s_1 + \lambda s_2)$ and $Z(\lambda) = \zeta'_{A,\Lambda}((1 - \lambda)s_1 + \lambda s_2)$, $0 \leq \lambda \leq 1$. Differentiating the first equation with respect to $\lambda$ we have: $z'(\lambda) = \zeta'_{A,\Lambda}((1 - \lambda)s_1 + \lambda s_2)(s_2 - s_1) = Z(\lambda)(1 - 2\sigma)$. We can always suppose that $1 - 2\sigma \geq 0$, in other words $s_1$ is at the left of $s_2$, since otherwise we can switch them.

If $1 - 2\sigma > 0$, we have $\arg z'(\lambda) = \arg Z(\lambda)$ , as long as $\zeta'_{A,\Lambda}((1 - \lambda)s_1 + \lambda s_2) \neq 0$, which means that the tangent to $\gamma$ at $z(\lambda)$ has the slope equal to that of the position vector of $Z(\lambda)$.

If $1 - 2\sigma = 0$, then $s_1 = s_2 = 1/2 + it$ and there is nothing to prove.

Having in view the color alternating rule, as well as the color matching rule, there are a few hypothetical positions of the curves $\Gamma_{k,j}$ and $\Upsilon_{k,j}$ which could allow such a configuration of zeros. We will show that every one of these hypothetical positions brings us to contradictions and therefore they need to be excluded. The simplest case is that of *embraced* curves $\Gamma_{k,j}$ which can appear for big values of $t$ (see [8]). In such a case the pre-image of the segment from $z = 1$ to $z = \zeta_{A,\Lambda}(v_{k,j})$ has as component an arc connecting $u_{k,j}$ and $u_{k,j+1}$ and passing through $v_{k,j+1}$ as shown in the second figure below. Due to the color alternating rule, it necessarily crosses once or an odd number of times the segment $I$. However, their images by $\zeta_{A,\Lambda}(s)$ are either disjoint or are intersecting an even number of times, which is a contradiction. A similar contradiction is obtained if we suppose that one of the ends of $I$ is the zero situated on $\Gamma_{k,0}$.

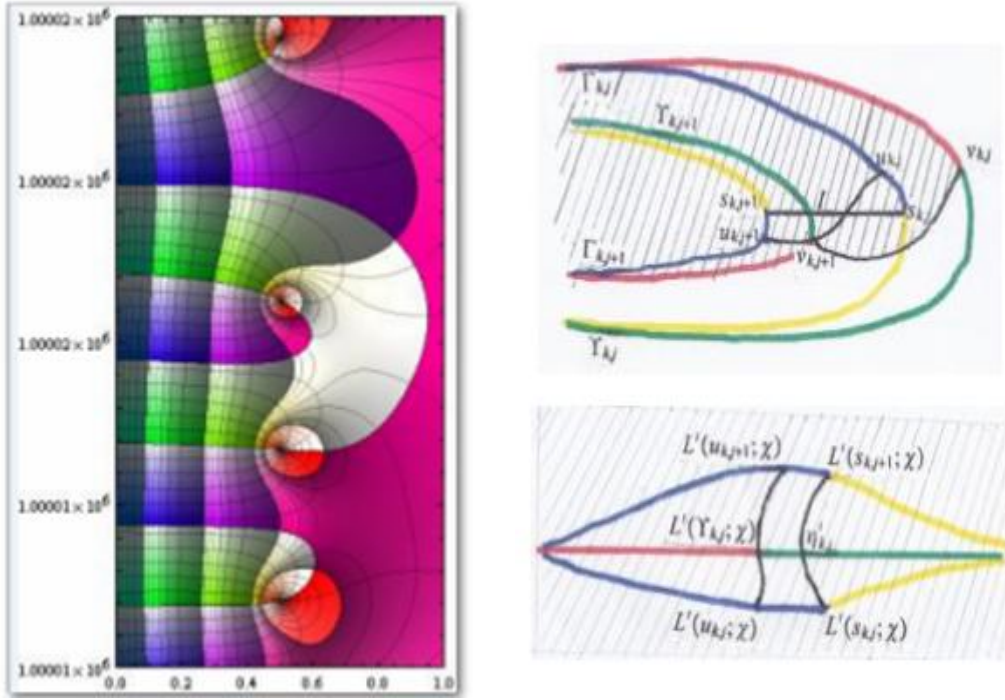

Fig.7

Therefore we only need to study the case where the ends of $I$ are the zeros $s_{k,j}$ at the left and $s_{k,j-1}$, respectively $s_{k,j+1}$ at the right, both above $\Gamma_{k,0}$ in the same strip $S_k$. Those zeros are separated by $\Upsilon_{k,j-1}$ and $\Upsilon_{k,j}$, respectively $\Upsilon_{k,j}$ and $\Upsilon_{k,j+1}$ as seen in Figs. 8 and 9 below. The situation where the zeros are below $\Gamma_{k,0}$ is symmetric and can be omitted. The part of $\Omega_{k,j}$ exterior to $\Gamma_{k,j}$ is mapped conformally by $\zeta_{A,\Lambda}(s)$ onto the upper half plane and the same is true for the part of $\Omega_{k,j-1}$ exterior to $\Gamma_{k,j-1}$. This means that $\gamma$ is situated in the upper half plane and its half-tangents at $0$ are one in the upper half plane and the other in the lower half plane. On the other hand, $\zeta'_{A,\Lambda}(s_{k,j})$ and $\zeta'_{A,\Lambda}(s_{k,j-1})$ are both in the upper half plane, which violates the equality $\arg z'(\lambda) = \arg Z(\lambda)$ for some values of $\lambda$, therefore such a configuration is impossible.

A different type of argument is necessary for the Fig. 9 below. At the points where $\arg z'(\lambda) = 0$ or $\arg z'(\lambda) = \pi$ the tangent to $\gamma$ is horizontal and $\gamma'$ crosses the positive, respectively the negative real half axis and vice-versa. We assigned numbers to those points. After $2$ the segment $I$ must cross $\Upsilon_{k,j+1}$ at $3$, which means that $\gamma'$ crosses the negative half axis at $3$ and correspondingly $\arg z'(\lambda)$ should be $\pi$. In the first case $\gamma$ crosses the positive real half axis at $4$ before going back to $1$, but then $\arg z'(\lambda)$ should take the value $\pi$, forcing $\gamma'$ to cross again the negative half axis, which contradicts the fact that $Z(1)$ should be in the upper half plane.

If the segment $I$ ends directly at $s_{k+1}$, then $\gamma$ goes back to $1$ remaining in the upper half plane and the half tangentat the end of $\gamma$ points towards the lower half plane.

Yet $0 < \arg Z(1) < \pi$.

Consequently, none of the hypothetical positions of the curves $\Gamma_{k,j}$ containing the two zeros is admissible and therefore $\sigma_1 = \sigma_2 = 1/2$. This shows that all the non trivial zeros of $\zeta_{A,\Lambda}(s)$ from $S_k$ must be situated on the line $\operatorname{Re} s = 1/2$. Since $S_k$ was arbitrary this is true for any non trivial zero.

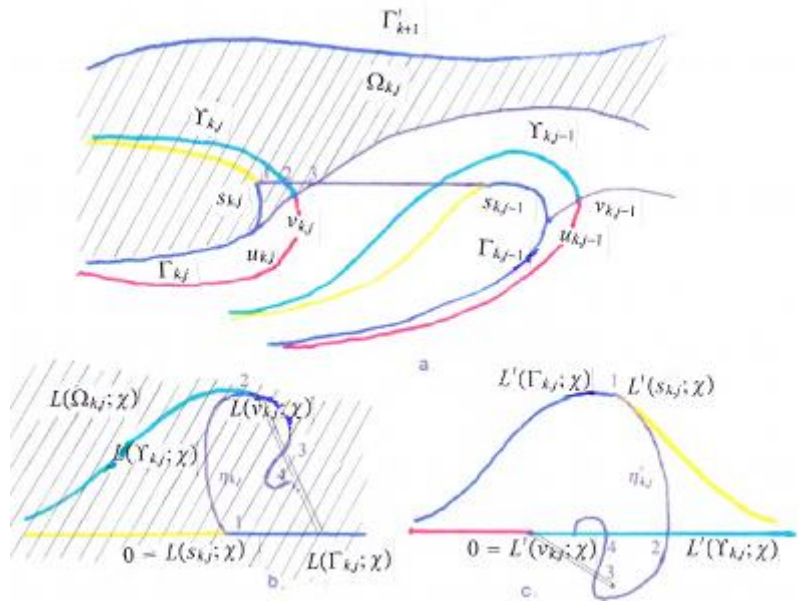

Fig. 8.

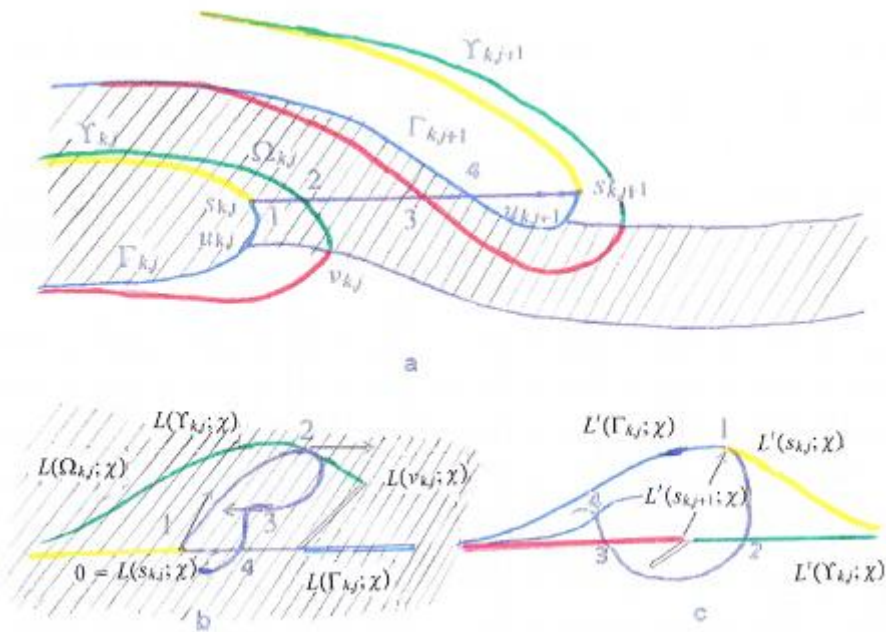

Fig. 9

## 5. The Bohr Functions Associated with a Dirichlet Series

Let $B = (\beta_n)$ be a basis for $\Lambda$ (see [3], section 8.3) and let $R$ be the Bohr matrix such that $\Lambda = RB$. By definition, the Bohr function corresponding to $B$ associated with the Dirichlet series (1) is

$$(11) \qquad F_B(Z) = \sum_{n=1}^{\infty} a_n Exp\{-(RZ)_n\},$$

where $(RZ)_n$ represents the n-th entry in the column matrix $RZ$ and $Z = (z_1, z_2, \ldots)$. We notice that if $z_n = s\beta_n$ for all $n$ i.e. $Z = sB$, then $RZ = sRB = s\Lambda$, hence $(RZ)_n = s\lambda_n$ and consequently

$$(12) \qquad F_B(sB) = \sum_{n=1}^{\infty} a_n e^{-s\lambda_n} = \zeta_{A,\Lambda}(s)$$

Then any zero $s_0$ of $\zeta_{A,\Lambda}(s)$ induces a zero $s_0 B$ of $F_B(Z)$. If $\zeta_{A,\Lambda}(s)$ verifies the hypothesis of Theorem 8, then we can localize the zeros of $F_B(Z)$ induced by the nontrivial zeros of $\zeta_{A,\Lambda}(s)$, namely:

**Theorem 9**. *All the zeros of $F_B(Z)$ induced by the non trivial zeros of $\zeta_{A,\Delta}(s)$ have every k-coordinate situated on the line* $\operatorname{Re} z_k = \frac{1}{2} b_k$ .

*Proof:* Indeed, if $s_0$ is a non trivial zero of $\zeta_{A,\Lambda}(s)$, then by Theorem 8, $\operatorname{Re} s_0 = 1/2$. Let $Z_0 = (z_1^{(0)}, z_2^{(0)}, \ldots)$ be the zero of $F_B(Z)$ induced by $s_0$. Then $\operatorname{Re} z_k^{(0)} = \operatorname{Re}(Bs_0)_k = \operatorname{Re}(b_k s_0) = b_k \operatorname{Re} s_0 = \frac{1}{2} b_k$.

[1] dghisa@yorku.ca

[2] lesandad.ferry@btinternet.com